\documentclass{amsart}

\usepackage{epsfig}
\usepackage{amsmath}
\usepackage{amssymb}
\usepackage{amscd}
\usepackage{graphicx}
\usepackage{color}
\usepackage{float}
\usepackage{verbatim}
\usepackage{bbm}
\usepackage{enumerate}

\newtheorem{thm}{Theorem}[section]
\newtheorem{lem}[thm]{Lemma}
\newtheorem{prop}[thm]{Proposition}

\newtheorem{cor}[thm]{Corollary}
\newtheorem{defn}[thm]{Definition}

\theoremstyle{remark}

\newtheorem{exam}[thm]{Example}

\def \N {\mathbb N}

\def \C {\mathfrak C}
\def \X {\mathfrak X}

\def \Z {\mathbb Z}
\def \R {\mathbb R}
\def \K {\mathcal K}
\def \F {\mathcal F}

\def \P {\mathcal P}

\def \M {\mathcal M}
\def \mtx {\mathcal M_T(X)}
\def \msy {\mathcal M_S(Y)}

\def \mtxe {\mathcal M^{\mathsf e}_T(X)}
\def \msx {\mathcal M_\sigma(X)}
\def \msxe {\mathcal M^{\mathsf e}_\sigma(X)}

\def \msup {\mathcal M_\sigma(\bar\X)}

\def \xt {$(X,T)$}
\def \xmt {$(X,\mu,T)$}

\def \xs {$(X,\sigma)$}

\def \us {$(\X,\sigma)$}
\def \ums {$(\X,\mu,\sigma)$}
\def \usp {$(\bar\X,\sigma)$}
\def \xtp {$(\bar X,\bar T)$}
\def \mtxp {\mathcal M_{\bar T}(\bar X)}
\def \mtxep {\mathcal M^{\mathsf e}_{\bar T}(\bar X)}
\def \eps {\varepsilon}
\def \fr {\mathsf{Fr}}
\def \supp {\mathsf{supp}}
\def \conv {\mathsf{conv}}

\def \per {\mathsf{Per}}
\def \aper {\mathsf{Aper}}
\def \uni {$\mathfrak X$}
\def \msu {\M_{\sigma}(\mathfrak X)}
\def \msue {\M^{\mathsf e}_{\sigma}(\mathfrak X)}

\def \sq {sequence}
\def \xsm {$(X,\Sigma,\mu)$}
\def \xsmt {$(X,\Sigma,\mu,T)$}

\def \xsmtp {$(\bar X,\bar\Sigma, \bar\mu,\bar T)$}

\def \ys {$(Y,S)$}
\def \tl {topological}
\def \im {invariant measure}
\def \inv {invariant}
\def \mpt {measure-preserving system}
\def \ds {dynamical system}
\def \zd {zero-dimensional}

\def \diam {\mathsf{diam}}
\def \Per {\mathsf{Per}}
\def \usc {upper semicontinuous}

\def \x {\mathfrak c}
\numberwithin{equation}{section}

\begin{document}

\title[Pure strictly uniform models of non-ergodic automorphisms]{Pure strictly uniform models of non-ergodic measure automorphisms}

\author{Tomasz Downarowicz and Benjamin Weiss}

\address{\vskip 2pt \hskip -12pt Tomasz Downarowicz}

\address{\hskip -12pt Faculty of Pure and Applied Mathematics, Wroc\l aw University of Technology, Wroc\l aw, Poland}

\email{downar@pwr.edu.pl}

\medskip
\address{\vskip 2pt \hskip -12pt Benjamin Weiss}

\address{\hskip -12pt Einstein Institute of Mathematics,
The Hebrew University of Jerusalem}

\email{weiss@math.huji.ac.il}

\thanks{
The first-named author is supported by
National Science Center, Poland (Grant HARMONIA No. 2018/30/M/ST1/00061) and
by the Wroc\l aw University of Science and Technology.
}


\begin{abstract}
The classical theorem of Jewett and Krieger gives a strictly ergodic model for any ergodic measure preserving system. An extension of this result for non-ergodic systems was given many years ago by George Hansel. He constructed, for any measure preserving system, a strictly uniform model, i.e. a compact space which admits an upper semicontinuous decomposition into strictly ergodic models of the ergodic components of the measure. In this note we give a new proof of a stronger result by adding the condition of purity, which controls the set of ergodic measures that appear in the strictly uniform model. 
\end{abstract}

\maketitle

\section{Introduction}
When a compact topological system \xt\ has a unique invariant measure $\mu$ it is easy to see that the measure preserving system \xsmt, where $\Sigma$ is the $\sigma$-algebra of Borel sets, is ergodic. Such a system is called \emph{uniquely ergodic}. It is also easy to see that when we restrict to the closed support $M$ of $\mu$, the resulting topological system is (in addition to uniquely ergodic) minimal, and such systems are called \textit{strictly ergodic}. More than 50 years ago R.\ I.\ Jewett \cite{J} proved the surprising theorem that any weakly mixing measure preserving system has a strictly ergodic model. Less than a year later  W.\ Krieger \cite{Kr} removed the extra condition of weak mixing and proved that every ergodic system has a strictly ergodic model. In 1974, G. Hansel \cite{H}  established a version of this result for non-ergodic systems.

Any non-ergodic system \xsmt\ has an ergodic decomposition. This is a measure $\xi$ on the space of ergodic measures on $X$, such that $\mu$ is the average with respect to $\xi$ of the ergodic measures in its support. Hansel showed that one can find a topological system \xt\ which is \emph{strictly uniform} and models the original system \xsmt. Strictly uniform means that the space $X$ decomposes as a union of strictly ergodic sets, moreover, the Ces\`aro means of every continuous function converge uniformly. Note that, by its nature, the support of $\xi$, as a subset of the set of ergodic measures for \xt, is only defined up to a null set, and \emph{a priori} the system \xt\ may carry ergodic measures completely unrelated to $\mu$ (however, $\xi$ will give measure zero to the set of such ``strange'' measures).

Our purpose in this note is to present a proof of a stronger version of the result of Hansel. We will show that given some set $\K$ of full measure with respect to $\xi$ one can arrange that the ergodic measures of the strictly uniform model belong~to~$\K$. In other words we can arrange that, for any choice for the support of $\xi$, only ergodic measures in that support will be found in the model. In Section \ref{sec2} we will give a precise formulation of our result and collect the basic tools that we will need. The construction itself is given in Section \ref{sec3}. In an earlier note \cite{DW} we gave a thorough discussion of strictly uniform systems and related classes of topological systems. In conclusion we would like to thank Mariusz Lema\'nczyk for raising a question which triggered this research.

\section{Pure strictly ergodic models; terminology and useful facts}\label{sec2}
\subsection{Topological models of \mpt s}

Throughout, by a \emph{\mpt} we will understand a quadruple \xsmt, where \xsm\ is a standard (Lebesgue) probability space and $T:X\to X$ is a measure automorphism. A \tl\ \ds\ \xt\ is a compact metric space $X$ with the action of a self-homeomorphism $T$. The set of $T$-\im s on $X$ will be denoted by $\mtx$ (always considered with the compact weak-star topology), while $\mtxe$ stands for the set of ergodic measures on $X$ (equivalently, extreme points of $\mtx$). Considered with a fixed \im\ $\mu\in\mtx$, \xt\ becomes a \mpt\ \xmt. The indication of the sigma-algebra is omitted intentionally, as this role will always be played by the sigma-algebra of Borel sets of $X$ (completed with respect to $\mu$, to create a standard probability space).

We will be using the following notation: If $\phi:(X,\Sigma)\to (\bar X,\bar\Sigma)$ is a measurable map between measurable spaces then $\phi^*$ denotes the adjoint map from the set $\M(X)$ of all probability measures on $(X,\Sigma)$ to $\M(\bar X)$ (defined analogously for $(\bar X,\bar\Sigma)$), given by
$$
\phi^*(\mu)(B)=\mu(\phi^{-1}(B)), \ \ (B\in\bar\Sigma).
$$

Two \mpt s, say \xsmt\ and \xsmtp, are said to be \emph{isomorphic} if there exists a measurable bijection $\phi:X'\to\bar X'$ between full subsets\,\footnote{A subset of a probability space is \emph{full} if it has measure 1.} of $X$ and $X'$, respectively, satisfying $\phi^*(\mu)=\bar\mu$, and which intertwines the actions of $T$ and $\bar T$, i.e.\ for which $\phi\circ T=\bar T\circ\phi$ $\mu$-almost surely. In such case, we will briefly say that $\mu$ and $\bar\mu$ are isomorphic.

A \tl\ \ds\ \xtp\ is called a \emph{model} for a \mpt\ \xsmt, if there exists an \im\ $\bar\mu\in\mtxp$ isomorphic to~$\mu$. If \xtp\ is also uniquely (respectively, strictly) ergodic, we call it a \emph{uniquely ergodic} (respectively, \emph{strictly ergodic}) \emph{model} of \xsmt. As we have already noted, if \xtp\ is a uniquely ergodic model of \xsmt, then $(\bar M,\bar T)$ is a strictly ergodic model of \xsmt, where $\bar M$ is a unique minimal subset of $\bar X$. We remark, that non-ergodic systems do not have uniquely ergodic models; they must carry enough ergodic measures to rebuild $\mu$.

A \tl\ \ds\ \xt\ is called \emph{uniform} if, for any continuous function $f:X\to\R$, the Ces\`aro means
$$
A_n^f=\frac1n\sum_{i=0}^{n-1}f\circ T^i
$$
converge uniformly on $X$. Uniform systems reveal a specific structure, captured in in the following statement, which is part of \cite[Theorem 4.9]{DW}:
\begin{thm}\label{doweis}
A \tl\ \ds\ \xt\ is uniform if and only if it admits an \usc\ partition whose atoms are closed, \inv\ and uniquely ergodic.
\end{thm}
If all atoms of the above mentioned partition are minimal (hence strictly ergodic), then the system is called \emph{strictly uniform}. G.\ Hansel (\cite{H}) proved that every \mpt\ \xsmt\ has a strictly uniform \zd\ model \xtp.

In this paper, we are going to give a new proof of Hansel's theorem, moreover, we are going to build a model with an additional property, which we call ``purity''. Roughly speaking, a model is pure if its collection of \im s is free of any ``expendable'' measures. In order to rigorously define purity we need to recall the notion of the ergodic decomposition. For technical reasons, it will be convenient to assume that the \mpt\ is already given in form of a \tl\ model, that is, we assume that we are given a system \xmt, where $\mu\in\mtx$.
\begin{defn}
Let $\mu\in\mtx$. The ergodic decomposition of $\mu$ is a probability measure $\xi$ on $\mtx$, such that
\begin{itemize}
	\item $\xi(\mtxe)=1$, and
	\item for any bounded measurable function $f:X\to\R$ one has
	$$
	\int f\,d\mu=\int\left(\int f\,d\nu\right)d\xi.
	$$
\end{itemize}
\end{defn}
It is well known that the ergodic decomposition of $\mu$ exists and is unique.

\begin{defn}
Let $\mu\in\mtx$ and let $\xi$ be its ergodic decomposition. Fix a set $\K\subset\mtxe$
satisfying $\xi(\K)=1$. A \tl\ \ds\ \xtp\ is called a \emph{pure (wrt.\ $\K$) model of $\mu$} if there exists a measurable bijection $\varphi:\mtxep\to\K'$, where $\K'\subset\K$ and $\xi(\K')=1$, and for every $\nu\in\mtxep$ the measures $\nu$ and $\varphi(\nu)$ are isomorphic.
\end{defn}

Observe that if \xtp\ is a pure model for $\mu$ then on $\mtxep$ we have the measure $\bar\xi=(\varphi^{-1})^*(\xi)$, and there exists a measure $\bar\mu\in\mtxp$, whose ergodic decomposition equals $\bar\xi$. Then $\bar\mu$ is isomorphic to $\mu$, which means that \xtp\ is indeed a \tl\ model of \xmt. It is convenient to think of the ergodic measures in $\K$ as ``admitted'' and those in the complement of $\K$ as ``unwanted''. A pure model supports the majority of the admitted measures without ``duplicating'', and no ``unwanted'' ergodic measures.

We can now formulate the main result of this paper.
\begin{thm}\label{main0}
With one exception, every \mpt\ \xsmt\ admits a strictly uniform \zd\ \tl\ model $(X_{\mathsf{psu}},\sigma)$, which is pure with respect to an a priori selected set $\K$ of\, $T$-ergodic measures on $(X,\Sigma)$, satisfying $\xi(\K)=1$, where $\xi$ is the ergodic decomposition of $\mu$.
\end{thm}
The above mentioned exception is a very special case of a purely periodic system described precisely in the ``working'' formulation of the main result, Theorem \ref{main}.
\medskip

Let us illustrate the concept of a pure model by three simple examples.

\begin{exam}
Let $X=[0,1]\times S^1$ (where $S^1=[0,1]/_{0=1}$) and let
$$
T(t,s)=(t,s+t\mod1).
$$ Let $\mu$ be the two-dimensional Lebesgue measure (which is clearly preserved by~$T$). The ergodic measures of the system \xt\ include irrational and rational rotations. The system is not aperiodic in the \tl\ sense\footnote{A \tl\ \ds\ is \emph{aperiodic} if it contains no fixed or periodic points.}, but the measure $\mu$ is aperiodic\footnote{An \im\ $\mu\in\mtx$ is \emph{aperiodic} if $\mu(\per)=0$, where $\per$ denotes the set of all $T$-periodic points (including fixpoints) in $X$.}; the ergodic decomposition of $\mu$ gives to the set of all rational rotations measure zero. So, we can choose $\K$ to be the set of all irrational rotations. A pure (wrt.\ $\K$) model of $\mu$ should have ergodic measures representing only the irrational rotations, each at most once. In particular, the pure model should be (in contrast to \xt) aperiodic in the \tl\ sense.
\end{exam}

\begin{exam}
Let $X=\{0,1\}^\Z$ and let $T$ be the shift transformation. Let $\mu$ be any non-ergodic \im\ with full topological support and zero entropy. Clearly, the ergodic decomposition of $\mu$ is supported by the set of all ergodic measures with entropy zero, so we can choose $\K$ to be this set. A pure model of \xmt\ will have (in contrast to \xt) \tl\ entropy zero.
\end{exam}

\begin{exam}
Let $X=\{0,1\}^\Z$ and let $T$ be the shift transformation. For $p\in(0,0.5]$, let $\nu_p$ be the Bernoulli measure with probabilities of $1$ and $0$ equal to $p$ and $1-p$, respectively. Let $\mu=\int_0^{0.5}\nu_p\,dp$. Fix $\K$ to be the set $\{\nu_p:p\in(0,0.5]\}$. A pure model of \xmt\ should have all ergodic measures isomorphic to Bernoulli measures, at most one for every entropy value $h\in(0,\log 2]$. In particular, it should have (in contrast to \xt) no ergodic measures of entropy zero.
\end{exam}

\subsection{Indistinguishable simplices of measures}
Given a metrizable simplex $\M$, a set $\F\subset\M$ is called a \emph{face} of $\M$ if $\F$ is itself a simplex, and all extreme points of $\F$ are extreme in $\M$. Observe that whenever $\K$ is a closed subset of the set of extreme points of $\M$, then the closed convex hull of~$\K$, $\overline{\conv(\K)}$, is a face of $\M$.\footnote{This need not be true if $\K$ is not closed. Although $\overline{\conv(\K)}$ is a closed convex subset of $\M$, it may have extreme points which are not extreme in $\M$ and it may fail to be a simplex.}

\begin{defn}
Let $\F_1$ and $\F_2$ be faces of $\mtx$ and $\mtxp$, for some \tl\ \ds s \xt\ and \xtp, respectively. We will say that $\F_1$ and $\F_2$ are \emph{indistinguishable} if there exists an affine homeomorphism $\varphi:\F_1\to\F_2$ such that the measures $\mu$ and $\varphi(\mu)$ are isomorphic, for every $\mu\in\F_1$.
\end{defn}

If $\F_1$ and $\F_2$ are Bauer simplices\footnote{A Bauer simplex is a simplex whose set of extreme points is compact.}, with $\K_1$ and $\K_2$ denoting the respective sets of extreme points, then $\F_1$ and $\F_2$ are indistinguishable if and only if there exists a homeomorphism $\varphi:\K_1\to\K_2$ such that $\mu$ and $\varphi(\mu)$ are isomorphic, for every $\mu\in\K_1$. In this case we will say that $\K_1$ and $\K_2$ are indistinguishable.\footnote{For indistinguishability of general (non-Bauer) simplices, indistinguishability of their sets of extreme points is insufficient. Take for instance the simplex of \im s of the subshift $X_1$ generated by the periodic \sq s \  $\dots0^{2k}1^{2k}0^{2k}1^{2k}\dots$ ($k\ge1$) and that of the subshift $X_2$ generated by the periodic \sq s\ $\dots0^{3k}1^k0^{3k}1^k\dots$ ($k\ge1$). Both simplices have two extreme points $\delta_0$ and $\delta_1$ supported by the fixpoints $\dots000\dots$ and $\dots111\dots$, and a \sq\ of periodic measures with periods $4k$ converging to a nonergodic measure. The sets of ergodic measures are homeomorphic and indistinguishable. However, in $X_1$ this \sq\ converges to $\frac{\delta_0+\delta_1}2$, while in $X_2$ the respective limit is $\frac{3\delta_0+\delta_1}4$, so the sets of \im s are not affinely homeomorphic.}

We will also use the notion of an isomorphic extension (which appears implicitly all over the literature, for an explicit appearance see \cite{DG}).

\begin{defn} We say that a \tl\ \ds\ \ys\ is an \emph{isomorphic extension} of \xt, if there exists a \tl\ factor map $\pi:Y\to X$ and a set $X'\subset X$ full for each \im\ on $X$, and such that $|\pi^{-1}(x)|=1$ for all $x\in X'$.
\end{defn}

Notice that in this case, the adjoint map $\pi^*:\msy\to\mtx$ is bijective (hence it is an affine homeomorphism) and for each $\mu\in\msy$, $\mu$ and $\pi^*(\mu)$ are isomorphic (the isomorphism is the map $\pi$). This implies that the simplices of measures $\mtx$ and $\msy$ are indistinguishable.

\subsection{The universal \zd\ system}\label{uzds}

Let $\X=\C^\Z$, where $\C$ is the Cantor set. If we equip $\X$ with the product topology,
it becomes a compact metrizable space homeomorphic to the Cantor set. The elements of $\X$ have the form of bilateral \sq s $\x=(\x_n)_{n\in\Z},\ \forall_{n\in\Z}\ \x_n\in\C$.
On $\X$ we consider the action of the \emph{shift transformation}, denoted by $\sigma$ and given by
$$
\sigma((\x_n)_{n\in\Z})=(\x_{n+1})_{n\in\Z}, \ \ (\x_n)_{n\in\Z}\in\X.
$$
The \tl\ \ds\ \us\ has several well-known \emph{universality properties}, which we list below:
\begin{enumerate}[(A)]
	\item Every \mpt\ \xsmt\ has a \zd\ \tl\ model, while every \zd\ \tl\ \ds\ \xt\ is
	conjugate to a subsystem (meaning a closed \inv\ subset) of the universal system.
  \item For every subsystem $(X,\sigma)$ of \us, $\msx$ is a face of $\msu$.
	\item \cite[Theorem~4.1]{D} For every face $\F$ of $\msu$ consisting of aperiodic
	measures there exists a closed \inv\ set $X\subset\X$ such that $\msx$ is
	indistinguishable from $\F$.
\end{enumerate}

\begin{figure}[H]
\includegraphics[width=13cm]{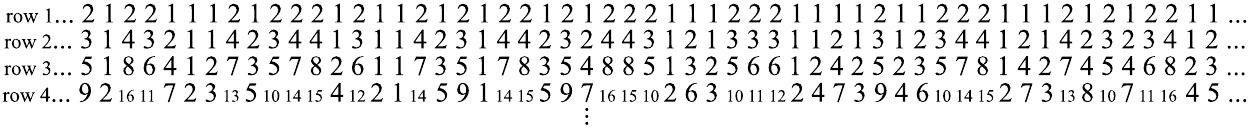}
\vspace{5pt}
\caption{An array $x\in\X$. Each symbol $x_{k,n}$ with $k\ge 2$ equals either $2x_{k-1,n}$ or $2x_{k-1,n}-1$.}
\label{fig0}
\end{figure}

\medskip
The universal system \uni\ has several conjugate representations. We will be using two of them. The first one, henceforth denoted by \us, is the \emph{inverse limit representation}:
$$
\mathfrak X = \overset\leftarrow{\lim_{k\to\infty}}\mathfrak X_k,
$$
where for each $k\ge1$, $\mathfrak X_k$ is the full shift over the alphabet $\Lambda_k=\{1,2,3,\dots,2^k\}$ and where the factor maps from $\mathfrak X_{k+1}$ onto $\mathfrak X_k$ are given by the ``amalgamations'' $\pi_k$ from $\Lambda_{k+1}$ to $\Lambda_k$, $\pi_k(m)=\lceil\frac m2\rceil$ (by abuse of notation, we will denote the resulting factor map by the same letter $\pi_k$ as the amalgamation). The system $\mathfrak X$ consists of all arrays
$$
x = [x_{k,n}]_{k\ge1,\,n\in\Z}
$$
satisfying, for every pair $(k,n)$, the relations $x_{k,n}\in\Lambda_k$ and
$x_{k,n}=\pi_k(x_{k+1,n})$. We think of $k$ as the row number and of $n$ as the column number. On our figures, the rows are ordered downward (that is, the first row appears on the top of an array). Figure \ref{fig0} shows an example of an array belonging to $\X$.

The transformation $\sigma$ acts on arrays as the horizontal left shift:
$$
\sigma([x_{k,n}]_{k\ge1,\,n\in\Z})=[x_{k,n+1}]_{k\ge1,\,n\in\Z}.
$$

\smallskip
The second representation, henceforth denoted by \usp, is the \emph{independent joining representation}:
$$
\bar\X = \prod_{k=1}^\infty\X_k,
$$
where, as before, for each $k\ge1$, $\X_k$ is the full shift over the alphabet $\Lambda_k$. This system also consists of all arrays
$$
x = [x_{k,n}]_{k\ge1,\,n\in\Z}
$$
satisfying, for every pair $(k,n)$, the relation $x_{k,n}\in\Lambda_k$, but this time no other restrictions are imposed. Observe that $\X$ is a subsystem of $\bar\X$, so whatever we say about the elements or \im s of $\bar\X$, it applies also to elements and \im s of $\X$.
\smallskip

By \emph{the $k$th row} of an element $x=[x_{k,n}]_{k\ge1,\,n\in\Z}\in\bar\X$ we will mean the $\Lambda_k$-valued sequence $x_{k,\cdot}=[x_{k,n}]_{n\in\Z}$.
We will also refer to the \sq\ $x_{\cdot,n}=[x_{k,n}]_{k\ge1}$ as the \emph{$n$th column} of $x$. If $x\in\bar\X$, $k\ge1$ and $m<n\in\Z$ then by
$x_{[1,k]\times[m,n]}$ we will denote the rectangular matrix
$$
x_{[1,k]\times[m,n]}=[x_{i,j}]_{1\le i\le k,\,m\le j\le n},
$$
and call it a \emph{rectangle (appearing in $x$) over $[1,k]\times[m,n]$}. Likewise, if $R$ is a rectangle over $[1,k]\times[m,n]$, and $1\le k'\le k$, $m\le m'<n'\le n$, then by $R_{[1,k']\times[m',n']}$ we will mean the ``sub-rectangle'' that occurs in $R$ over  $[1,k']\times[m',n']$. In what follows, equality between rectangles will always be understood modulo horizontal shift.
\smallskip

Any rectangle $R$ over $[1,k]\times[m,n]$ determines an \emph{empirical measure} on all rectangles $Q$, by frequencies, as follows: If $Q$ is over $[1,l]\times[s,t]$ and
either $l>k$ or $n-m>t-s$ then we let $\fr_R(Q) = 0$. Otherwise we let 
$$
\fr_R(Q) = \frac1{(n-m)-(t-s)}|\{i\in[m,n-(t-s)]:R_{[1,l]\times[i,i+(t-s)]}=Q\}|.
$$
Clearly, the frequency of $Q$ does not depend on the shift, so we can restrict the definition to rectangles $Q$ over $[1,l]\times[0,t-s]$ called briefly \emph{rectangles of dimensions} $l\times(t-s+1)$.
The frequencies allow to define a distance between rectangles $R$ and \im s $\mu\in\msup$, as follows:
$$
d^*(R,\mu) = \sum_{l\ge1, r\ge 1}2^{-l-r}\sum_{\text{rectangles $Q$ of dimensions $l\times r$}} |\fr_R(Q)-\mu([Q])|,
$$
where by $[Q]$ we mean the cylinder $\{x\in\bar\X: x_{[1,k]\times[0,r-1]}=Q\}$.
\smallskip

It follows from the ergodic theorem that if $\mu$ is an ergodic measure carried by a subsystem $X\subset\bar\X$ then it can be approximated by rectangles appearing in $X$. Conversely, given $\eps>0$, any rectangle of sufficiently large dimensions, appearing in $X$, is $\eps$-close to some \inv\ (however not necessarily ergodic) measure on $X$.

By an obvious modification, the metric $d^*$ establishes also a distance between two \im s as well as between two rectangles. It is clear that $d^*\le 1$ and $d^*$ is a convex metric on \im s. The distance between rectangles and measures satisfies an approximate version of convexity, which we formulate below (we skip the standard proof):

\begin{prop}\label{apcon} Suppose that $R_l$, $l=1,2,\dots,q$, $q\in\N$, are rectangles of dimensions $k\times n_l$ such that
$$
d^*(R_l,\mu_l)<\eps,
$$
for some \im s $\mu_l$. Let $R$ denote the concatenation $R_1R_2\dots R_q$ ($R$ is a rectangle of dimensions $k\times n$, where $n=n_1+n_2+\cdots+n_q$). Let $\bar\mu = \sum_{l=1}^q \frac {n_l}n \mu_l$. Then
$$
d^*(R,\bar\mu)<\eps+O(\tfrac qn).
$$
\end{prop}

\subsection{Markers in aperiodic systems}

By \emph{putting a marker} in an array $x=[x_{k,n}]_{k\ge1,\, n\in\Z}\in\bar\X$ at the position $(k,n)$ we will understand placing an additional symbol ``$|$'' in row $k$ between the symbols appearing at the coordinates $n$ and $n+1$. To allow markers, we formally need to replace each alphabet $\Lambda_k$ used in row $k$ by $\hat\Lambda_k = \Lambda_k\times\{\emptyset, | \}=\{1,2,\dots,2^k\}\cup\{1|,2|,\dots,2^k|\}$.

\begin{defn}\label{distma}We let $(\hat\X,\sigma)$ denote the extension of the universal system \usp\ obtained by putting, in every row of every array $x\in\bar\X$, markers, in all possible ways satisfying the conditions (1)-(4) below (in this manner creating multiple preimages $\hat x\in\hat\X$ of each $x\in\bar\X$):
\begin{enumerate}
	\item \emph{Shift equivariance}: If $\hat x$ has a marker at a position $(k,n)$, then
	$\sigma(\hat x)$ has a marker at the position $(k,n-1)$.
	\item \emph{Two gap sizes}: There exists a (quickly) increasing \sq\ of integers
	$(l_k)_{k\ge 1}$ such that, in every $\hat x\in\hat\X$, the markers in row $k$ appear with
	only two gap sizes, $l_k$ and $l_k+1$.
	\item \emph{Balanced frequencies of gaps}: There exist integers $L_k$, $k\ge1$, such that,
	for each $\hat x$, in each interval of length $L_k$ in row $k$ of $\hat x$, there appear
	at least $\frac{L_k}{3l_k}$ gaps of length $l_k$ and at least $\frac{L_k}{3(l_k+1)}$ gaps
	of length $l_k+1$.
	\item \emph{Congruency}: For any $\hat x\in\hat\X$, the set of positions of markers in row $k\!+\!1$ is a subset of the set of positions of markers in row $k$.
\end{enumerate}	
\end{defn}

\begin{figure}[H]
\includegraphics[width=13cm]{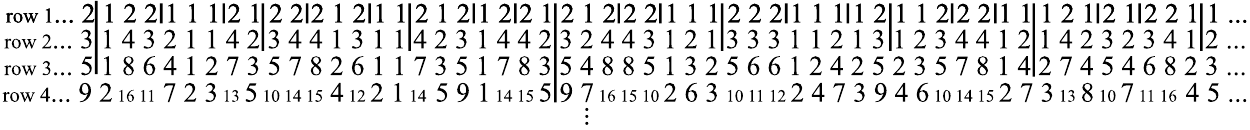}
\vspace{5pt}
\caption{An array $\hat x\in\hat\X$.}
\label{fig01}
\end{figure}

In the rest of this paper all subsystems of $\hat\X$ will be denoted with a ``hat'', e.g.~$\hat X$. Figure \ref{fig01} shows an array (the same as that on Figure \ref{fig0}) equipped with a system of markers satisfying (2), (3) and (4) (condition (1) cannot be seen in one array).

\begin{defn}We say that a \tl\ \ds\ \xt\ is \emph{measure-saturated} if it equals the closure of the union of \tl\ supports of all its \im s (this closure is often called the \emph{measure center}).
\end{defn}

The following lemma is a \tl\ analog of what is known in ergodic theory as the Kakutani--Rokhlin Lemma (specific version with only two tower heights).

\begin{lem}\label{marker2}
Every aperiodic, zero-dimensional, measure-saturated system \xs, where $X\subset\bar\X$, admits an isomorphic extension $(\hat X,\sigma)$ with $\hat X\subset\hat\X$, where the factor map $\pi:\hat X\to X$ consists in erasing the markers.
\end{lem}

\begin{proof}Because the proof is standard, we only outline its main ideas skipping some technical details.

Our basic tool is Krieger's marker lemma for \zd\ systems (see~\cite[Lemma~2.2]{B}) which, applied to aperiodic systems, implies that for any natural $N>1$ there exists a clopen set (called the \emph{$N$-marker set}), visited by every orbit with gaps ranging between $N$ and $2N-1$. We can choose a fast growing \sq\ $(N_k)_{k\ge1}$ and in every array $x\in X$ place markers in the $k$th row at the positions corresponding to the visit times of the orbit of $x$ in the $N_k$-marker set. Because the marker sets are clopen, this process is a \tl\ conjugacy. The temporary ``markered'' array system $(\check X,\sigma)$ that we get satisfies the following conditions:
\begin{enumerate}[(a)]
	\item \emph{Shift equivariance}: If $\check x\in \check X$ has a marker at a position $(k,n)$, then
	$\sigma(\check x)$ has a marker at the position $(k,n-1)$.
	\item\emph{Controlled gap sizes}: The markers in row $k$ appear with gaps ranging between
	$N_k$ and $2N_k-1$.
	\item The map $\check x\to x$ consisting in erasing all markers is a conjugacy between
	$\check X$~and~$X$.
\end{enumerate}	
We continue as follows. Given $l_1$, we can assume that $N_1$ is larger than $9l_1^2+l_1$. For each number $p\ge N_1$ we fix a pair of positive integers $a_p$ and $b_p$ such that  $\frac89<\frac{a_p}{b_p}<\frac98$ and $p=a_pl_1+b_p(l_1+1)$ (it is an easy exercise to show that such integers $a_p$ and $b_p$ exist). Now, in row number 1 we introduce additional markers which subdivide each gap between the original markers into $a_p$ pieces of length $l_1$ followed by $b_p$ pieces of length $l_1+1$, where $p$ is the length of the gap.
The resulting system $\check X_{(1)}$ is conjugate to $X$. Indeed, we already know that $x\mapsto \check x$ is a conjugacy. Further, the algorithm of introducing the additional markers is deterministic, continuous and shift-equivariant, and thus it is a factor map $\check x\mapsto\check x_{(1)}\in \check X_{(1)}$. Finally, erasure of all markers, $\check x_{(1)}\mapsto x$, is a factor map from $\check X_{(1)}$ to $X$, inverse to the composition $x\mapsto\check x\mapsto\check x_{(1)}$.
The inequalities $\frac89<\frac{a_p}{b_p}<\frac98$ guarantee that in every interval of length $L_1=6N_1$ at least one third of the gaps have length $l_k$ and at least one third have length $l_k+1$.

We pass to step 2, assuming that $l_2$ is much larger than $l_1$ and that $N_2$ is larger than $9l_2^2+l_2$. First, in all elements $\check x_{(1)}\in\check X_{(1)}$ we place additional markers in row~2, so that the gaps between them are either or $l_2$ and $l_2+1$. The algorithm follows identical rules as in step~1. However, the markers in row 2 need not satisfy the congruency condition with respect to the markers in row 1. To fix this, we need to revise the markers in row~1 (maintaining the gaps lengths $l_1$ and $l_1+1$ and their proportion in every interval of length $L_1$). To achieve this, it suffices to rearrange the markers within the distance $9l_1^2+l_1$ on both sides of each marker in row~2. Because $l_2$ is much larger than $l_1$, the modifications of row~1 affect only a very small percentage (say $\eps_1$ in terms of upper Banach density) of the markers in that row. Moreover, we can do that via a deterministic procedure, i.e.\ for each configuration of markers in row~1 within the distance $9l_1^2+l_1$ on both sides of each marker in row~2 (which misses a marker in row~1) we choose one arrangement of markers in row~1 that fixes the miss and apply it whenever that configuration is observed. In this manner we have constructed a system $\check X_{(2)}$, which, for identical reasons as in step 1, is conjugate to $X$.

In step 3, in all elements $\check x_{(2)}$ of $\check X_{(2)}$, we insert new markers in row 3, so that the new gaps are $l_3$ and $l_{3+1}$, using the same rule as in steps 1 and 2, and next we rearrange a small percentage (say $\eps_2$) of markers in row 2 to fix the congruency between markers in rows 3 and 2. Then we need to rearrange markers in row 1 again, to fix the congruency between markers in rows 2 and 1, but this affects only the percentage $\eps_2$ of markers in row~1. We believe that further inductive construction is now understood. We need to ensure the growth of the numbers $l_k$ so fast that the resulting \sq\ $(\eps_k)_{k\ge1}$ is summable.

By the above mentioned summability and a topological version of the Borel--Cantelli Lemma (in which measure is replaced by the upper Banach density), there is a set $X_0\subset X$ full for each \im\ on $X$, such that for every $ x\in X_0$ and each row $k$, every marker established in that row in the inductive step $k$, is affected (moved) at most finitely many times throughout the rest of the construction. So, for each point $x\in X_0$, the systems of markers introduced throughout the induction converge to a system of markers in all rows that obeys the conditions (1)--(4) of Definition \ref{distma}. We let $\hat x=\lim\check x_{(k)}$ (coordinatewise limit) denote the point $x\in X_0$ equipped with the above system of markers. In our final move, we let $\hat X$ be the closure of the set of arrays $\{\hat x:\ x\in X_0\}$. Obviously, the properties (1)--(4) pass to all elements of $\hat X$, hence $\hat X\subset\hat\X$. Note that the erasure of all markers is a \tl\ factor map from $\hat X$ onto the closure of $X_0$, which, by the measure-saturation of $X$, equals $X$. It is not hard to verify that points in $X_0$ have a unique preimage, namely the array $\hat x$ constructed inductively as described above, and only the points from $X\setminus X_0$ may obtain in $\hat X$ multiple configurations of markers. So, $(\hat X,\sigma)$ is an isomorphic extension of \xs. This ends the proof.
\end{proof}

\subsection{Approximation of measures by $k$-rectangles}
\begin{figure}[H]
\includegraphics[width=13cm]{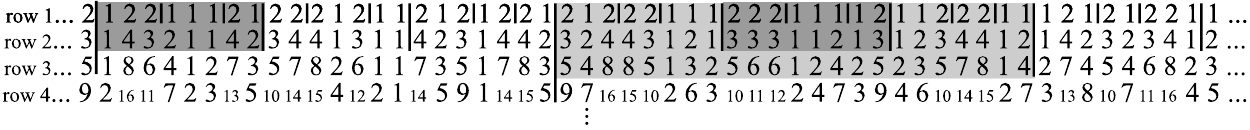}
\vspace{5pt}
\caption{Selected $k$-rectangles from the array on Figure~\ref{fig01} (two $2$-rectangles shaded dark-gray, and one $3$-rectangle shaded light-gray).}
\label{fig1}
\end{figure}

Let $(\hat X,\sigma)$ be the ``markered'' isomorphic extension of an aperiodic \zd\ measure-saturated system \xs, constructed in Lemma \ref{marker2}.
\begin{defn}
By a $k$-rectangle we will mean the rectangular array (with markers) of dimensions $k\times l_k$ or $k\times (l_k\!+\!1)$, appearing in rows $1$ thorugh $k$ of some $\hat x\in\hat X$, between two consecutive markers in row $k$ (see Figure \ref{fig1}).
\end{defn}

\begin{defn} Fix a positive number $\gamma$. Let $\hat\mu\in\M^{\mathsf e}_\sigma(\hat X)$ be an ergodic measure on $\hat X$. A $k$-rectangle $R$ is \emph{$\gamma$-close to $\hat\mu$} if
$$
d^*(R,\hat\mu)<\gamma
$$
(the presence of markers does not affect the definition of $d^*$; we only operate on larger yet finite alphabets).
\end{defn}

\begin{lem}\label{goodrec} Consider an ergodic measure $\hat\mu\in\M^{\mathsf e}_\sigma(\hat X)$ and  positive numbers $\gamma$ and $\eps$. The following facts hold:
\begin{enumerate}
\item Let $A_{k,\hat\mu,\gamma}$ denote the set of points $\hat x\in\hat X$, such that the \emph{central} (i.e., covering the coordinate zero) $k$-rectangle in $\hat x$ is $\gamma$-close to $\hat\mu$. Then, for any sufficiently large integer $k\in\N$, we have $\hat\mu(A_{k,\hat\mu,\gamma})>1-\eps$.
\item If $\hat\mu(A_{k',\hat\mu,\frac\gamma2})>1-\frac{\gamma^2}2$ for large enough $k'$ then
\begin{enumerate}
	\item for any $k>k'$, we have $\hat\mu(A_{k,\hat\mu,\gamma})>1-\gamma$,
	\item there are $k$-rectangles of both lengths $l_k$ and $l_k\!+\!1$ which are
	$\gamma$-close~to~$\hat\mu$.
\end{enumerate}
\end{enumerate}
\end{lem}

\begin{proof}
Given integers $k_0,m_0$, let $A_{k_0,m_0}$ denote the set of points $\hat x\in\hat X$ satisfying the following condition: for any $m\ge m_0$ the rectangles
$\hat x_{[1,k_0]\times[0,m-1]}$ and $\hat x_{[1,k_0]\times[-m+1,0]}$ (of dimensions $k_0\times m$) are $\frac\gamma3$-close to $\hat\mu$.

By the ergodic theorem, there exists a set $A\subset \hat X$ with $\hat\mu(A)=1$ of points (i.e.\ arrays) generic for $\hat\mu$ under both $\sigma$ and $\sigma^{-1}$. This easily implies that there exist (arbitrarily large) integers $k_0$ and $m_0$ such that $\hat\mu(A_{k_0,m_0})>1-\eps$. Because rows with large indices $k$ have small influence on the metric $d^*$, if $k_0$ is large enough, then for any $k\ge k_0$ (and $m\ge m_0$) the rectangles $\hat x_{[1,k]\times[0,m-1]}$ and $\hat x_{[1,k]\times[-m+1,0]}$, where $\hat x\in A_{k_0,m_0}$, are $\frac{2\gamma}3$-close to $\hat\mu$. Let $k\ge k_0$ be so large that $l_k$ is much larger than $m_0$. Consider the central $k$-rectangle $R$ in an array $\hat x\in A_{k_0,m_0}$. There are three possibilities:
\begin{itemize}
	\item $R$ extends far enough to the left and right so it covers (in the horizontal
	direction) the coordinates $[-m_0,m_0]$. Then both its parts (left and right from
	coordinate zero) are $\frac{2\gamma}3$-close to $\hat\mu$, and then $R$ is $\gamma$-close to
	$\hat\mu$.
	\item $R$ starts at a coordinate $n\in(-m_0,0)$. Then the right part of $R$ is
	$\frac{2\gamma}3$-close to $\hat\mu$, while its left part is negligible in comparison to the
	right part. In this case	$R$ is also $\gamma$-close to $\hat\mu$
	\item $R$ ends at a coordinate $n\in(0,m_0)$. Then the left part of $R$ is
	$\frac{2\gamma}3$-close to $\hat\mu$, while its right part is negligible and the conclusion is
	the same as above.
\end{itemize}
We have shown that the set of arrays whose central $k$-rectangle is $\gamma$-close to $\hat\mu$ contains $A_{k_0,m_0}$, which ends the proof of (1).

Now assume that for some large $k'$ we have $\hat\mu(A_{k',\hat\mu,\frac\gamma2})>1-\frac{\gamma^2}2$, and fix some $k>k'$. Any $k$-rectangle is a concatenation of $k'$-rectangles (and some contents added in rows $k'\!+\!1,k'\!+\!2,\dots,k$; this contents can be ignored when calculating distances). Let us say that a $k$-rectangle $R$ is ``OK'', if the $k'$-rectangles which are $\frac\gamma2$-close to $\hat\mu$ contribute more than the fraction $1-\frac\gamma2$ to the length of $R$. Proposition~\ref{apcon} implies that $d^*(R,\hat\mu)$ may exceed the convex combination $(1-\frac\gamma2)\cdot\frac\gamma2 + \frac\gamma2\cdot1$ (which is strictly smaller than $\gamma$) only by a small fraction, and thus, if $k'$ is large enough, we have $d^*(R,\hat\mu)<\gamma$. This means that whenever the central $k$-rectangle of some array $\hat x$ is ``OK'' then $\hat x\in A_{\hat\mu,k,\gamma}$.

The condition $\hat\mu(A_{k',\hat\mu,\frac\gamma2})>1-\frac{\gamma^2}2$ implies that in a $\hat\mu$-typical array $\hat x$ the density of the set of the horizontal coordinates occupied by the $k'$-rectangles which are $\frac\gamma2$-close to $\hat\mu$ exceeds $1-\frac{\gamma^2}2$. This easily implies that the density of the set of coordinates occupied by the $k$-rectangles which are not ``OK'' must not exceed $\gamma$. Translated to the terms of measure, we have shown that the set of arrays $\hat x$ whose central $k$-rectangle is ``OK'' exceeds $1-\gamma$ in measure. Since this set contains $A_{k,\hat\mu,\gamma}$, (2a) is proved.

For (2b), it suffices to choose $\gamma<\frac13$. By the condition~(4) in Definition \ref{distma} (balanced frequencies of gaps), for any $k$, the set of arrays $\hat x$ whose central $k$-rectangle has length $l_k$ has measure $\hat\mu$ between $\frac13$ and $\frac23$ (and the same is true for the length $l_k\!+\!1$). So, if $k$ is so large that~(2a) holds for $\gamma<\frac13$, both of these sets have nonempty intersection with the set $A_{k,\hat\mu,\gamma}$ of arrays whose central $k$-rectangle is $\gamma$-close to $\hat\mu$. Thus $k$-rectangles $\gamma$-close to $\hat\mu$ of both lengths exist.
\end{proof}

\begin{cor}\label{corol}
If $\hat\K\subset\M_\sigma(\hat X)$ is compact then for any $\gamma>0$ there exists $k$ such that for all $\hat\mu\in\hat\K$ the set $\hat\mu(A_{k,\hat\mu,\gamma})$ has measure $\hat\mu$ at least $1-\gamma$ and there exist $k$-rectangles of both lengths $l_k$ and $l_k\!+\!1$ which are $\gamma$-close to $\hat\mu$.
\end{cor}

\begin{proof}
By Lemma \ref{goodrec} (1), for each $\hat\mu\in\hat\K$ there exists (an arbitrarily large) $k'$ such that $\mu(A_{k',\hat\mu,\frac\gamma2})>1-\frac{\gamma^2}2$. Since the set $A_{k',\hat\mu,\frac\gamma2}$ is open and depends continuously on $\hat\mu$, this inequality (with a fixed parameter $k'$) holds on an open set of measures. By compactness, there are finitely many values of $k'$ which suffice for all measures $\hat\mu\in\hat\K$.
By Lemma \ref{goodrec} (2), any $k$ larger than the maximum of the above mentioned finitely many values of $k'$, satisfies the desired condition.
\end{proof}

\section{Pure strictly ergodic models; the construction}\label{sec3}

Let \xsmt\ be a \mpt\ which is not ergodic (due to Jewett--Krieger Theorem, ergodic systems are, from our point of view, trivial). According to Theorem \ref{main0}, we intend to built a \tl\ \zd\ model $(X_{\mathsf{psu}},\sigma)$ of \xsmt, which is strictly uniform and pure with respect to an a priori selected set $\K$ of admitted ergodic measures. By the ``universality property''~(A) of the universal system, we may assume that \xsmt\ is given as a \zd\ \tl\ model, moreover, that model is a subsystem of the universal system \us\ (the inverse limit representation). In fact, we can assume that the model equals \ums, where $\mu\in\msu$ and hence the ergodic decomposition $\xi$ of $\mu$ is supported by $\msue$, and $\K$ is a subset of $\msue$.

\subsection{Modeling a Cantor set of ergodic measures}

The pure strictly uniform model in the general case will be built in several stages.
The most important and technically crucial is the construction of a pure uniform (but not yet strictly uniform) model in case $\K$ is homeomorphic to the Cantor set and contains no periodic measures.

\begin{thm}\label{case1}
Let $\K\subset\msue$ be homeomorphic to the Cantor set, and contain no periodic measures. Then there exists a uniform \zd\ system $(\hat X_{\mathsf{pu}},\sigma)$ such that $\msxe$ is indistinguishable from $\K$.
\end{thm}

\begin{proof}
Since $\K$ is closed, $\F=\overline{\conv}(\K)$ is a face in $\msue$. By assumption, $\F$ consists of aperiodic measures only. By (C) in subsection \ref{uzds}, there exists a \zd\ \tl\ \ds\ \xs\ such that $\msx$ is indistinguishable from $\F$, which, in case of a Bauer simplex, is equivalent to $\msxe$ being indistinguishable from $\K$. Thus, we can as well assume that $\K=\msxe$. By universality of the inverse limit system \us, we may assume that $X\subset\X$. Restricting $X$ to its measure center, we can also assume that \xs\ is measure-saturated. Because $\X\subset\bar\X$, we can use Lemma \ref{marker2} which implies that \xs\ has an isomorphic extension $(\hat X,\sigma)$ contained in $\hat\X$ ($\hat X$ has markers satisfying the conditions (1)-(4) of Definition \ref{distma}). We let $\phi:\hat X\to X$ denote the corresponding factor map consisting in erasing the markers. Clearly, $\K$ and $\hat\K=\M^{\mathsf e}_\sigma(\hat X)$ are indistinguishable. We let
$\Phi=(\phi^*)^{-1}:\K\to\hat\K$ be the inverse of the map adjoint~to~$\phi$.
We establish a summable \sq\ of positive numbers $(\eps_m)_{m\ge 1}$. The construction of $\hat X_{\mathsf{pu}}$ goes by induction~on~$m$.
\smallskip

{\bf Step 1}. We let $\P_1=\{K_1,K_2,\dots,K_{r_1}\}$ be a clopen partition of $\K$ into sets so small that the their images by $\Phi$, denoted by $\hat K_i$ (which form a clopen partition of $\hat\K$) have diameters smaller than $\eps_1$. Since each set $\hat K_i$ consist of ergodic measures, the closed convex hulls of the sets $\hat K_i$ are disjoint.
We let $\gamma_1<\eps_1$ be smaller than half of the smallest distance between the above mentioned closed convex hulls. By Corollary~\ref{corol} and Proposition \ref{apcon},
there exists $k_1$ so large that
\begin{enumerate}
	\item[(a)] for each $\hat\mu\in\hat\K$ the measure of the set of arrays $\hat x$ whose
	central $k_1$-rectangle is $\gamma_1$-close to $\hat\mu$ exceeds $1-\gamma_1$,
	\item[(b)] there exists at least one $k_1$-rectangle of each length $l_{k_1}$ and $l_{k_1}\!+\!1$, which is $\gamma_1$-close to $\hat\mu$,
	\item[(c)] if $k_1$-rectangles $R_1,\dots R_q$ satisfy $d^*(R_l,\hat\mu_l)<\gamma_1$, for some \im s $\hat\mu_l$, $l=1,2,\dots,q$, $q\in\N$, then $d^*(R,\bar\mu)<\eps_1$, where $R$ denotes the concatenation $R_1R_2\dots R_q$ and $\bar\mu$ is some convex combination of the measures~$\hat\mu_l$.
\end{enumerate}

The further procedure will be described for a fixed index $i=1,2,\dots,r_1$. Each of the $k_1$-rectangles appearing in $\hat X$, falling in the $\gamma_1$-neighborhood of $\hat K_i$ will be called \emph{$i$-good}, all other $k_1$-rectangles will be called \emph{$i$-bad}. In other words, a $k_1$-rectangle is $i$-good if it is $\gamma_1$-close to some $\hat\mu\in\hat K_i$. By (b), there exist $i$-good $k_1$-rectangles of each length $l_{k_1}$ and $l_{k_1+1}$. We select one $i$-good $k_1$-rectangle of each length and we call these two rectangles the \emph{$i$-tabbed $k_1$-rectangles}. Next, in every array $\hat x\in\hat X$ we replace all $i$-bad $k_1$-rectangles by the $i$-tabbed $k_1$-rectangles of matching lengths
(see Figure \ref{fig3}).
\begin{figure}[H]
\includegraphics[width=13cm]{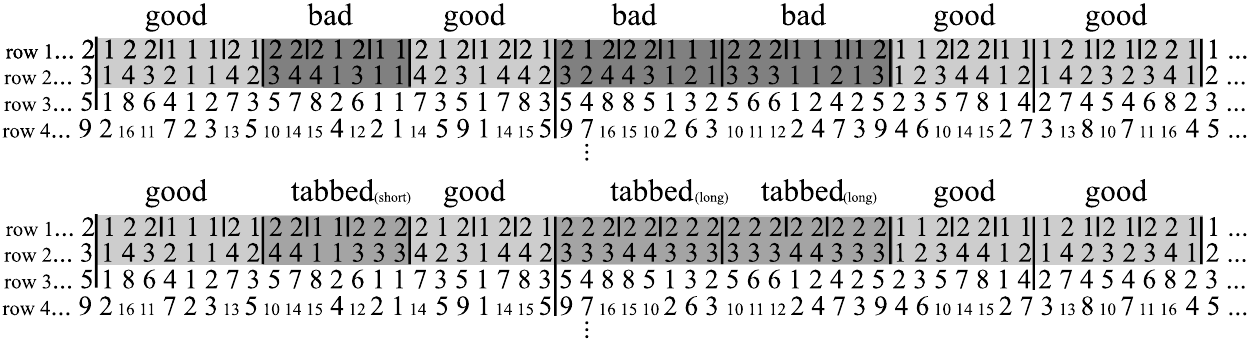}
\caption{The top figure shows the classification of $2$-rectangles into good and bad.
The bottom figure shows bad rectangles replaced by the tabbed rectangles of the same size. Note that some markers in row 1 have moved (but this movement does not affect the construction).}
\label{fig3}
\end{figure}
The modified system $\hat X$ will be denoted by $\hat X_i$ and the modification itself will be denoted by $\theta_i: \hat X\to \hat X_i$. By (a), the set of arrays whose central $k_1$-rectangle has been changed by $\theta_i$ has, for each $\hat\mu\in\hat K_i$, a measure smaller than~$\gamma_1$ (which is less than $\eps_1$).

At the same time, we denote by $\hat S_i$ the closed \inv\ subset of $\hat\X$ consisting of all possible arrays allowing any concatenations of $i$-good $k_1$-rectangles in the first $k_1$-rows (and anything allowable in other rows). It is clear that $\hat X_i\subset \hat S_i$. Since, by the choice of $\gamma_1$, the families of $i$-good $k$-rectangles are disjoint for different indices~$i$, the sets $\hat S_i$ are disjoint as well.

The map $\theta_i$, being a code with finite horizon, is continuous and shift-equivariant, so $\hat X_i$ is a \tl\ factor of $\hat X$. On the other hand, $\theta_i$ leaves the rows with indices higher than $k_1$ intact, therefore there exists a factor map $\phi_i$ from $\hat X_i$ onto the original system $X$ (without the markers). Indeed, $\phi_i$ consists in erasing all markers and reconstructing the original rows $1$ through $k$ using the (unchanged) $(k\!+\!1)$st row and the successive amalgamations. As easily verified,
$$
\phi_i\circ\theta_i=\phi,
$$
which implies that $\hat X_i$ is an isomorphic extension of $X$. The inverse of the adjoint map, $\Phi_i=(\phi_i^*)^{-1}:\K\to\hat\K_i=\M^{\mathsf e}_\sigma(\hat X_i)$ is an affine homeomorphism sending measures to their isomorphic copies. As the sets $\hat\K_i$ consist of measures supported by the disjoint systems $\hat S_i$, they are disjoint. Each of them is indistinguishable~from~$\K$.

Because all elements of $\hat S_i$ have in the first $k_1$ rows concatenations of only $i$-good $k_1$-rectangles, by (c), all \im s supported by $\hat S_i$ (in particular all elements of $\hat\K_i$) lie within the $\eps_1$-neighborhood of the closed convex hull of $\hat K_i$. By convexity of $d^*$ this hull itself has diameter smaller than $\eps_1$. Therefore, for each $\mu\in K_i$, we have
$$
d^*(\Phi(\mu),\Phi_i(\mu))<2\eps_1.
$$
\smallskip

{\bf Step $m\!+\!1$}. Assume that for some $m\ge1$ we have chosen integers $r_1,r_2,\dots,r_m$ and a partition $\P_m$ of $\K$ whose atoms are clopen sets not exceeding $\frac{\diam(\K)}m$ in diameter, denoted by  $K_{i_1,i_2,\dots,i_m}$, where $(i_1,i_2,\dots,i_m)\in\prod_{j=1}^m\{1,2,\dots,r_j\}$.
Moreover, if $m\ge 2$, we assume that $K_{i_1,i_2,\dots,i_m}\subset K_{i_1,i_2,\dots,i_{m-1}}$, for each string $(i_1,i_2,\dots,i_m)$.
Next, suppose that we have selected some integers $k_1<k_2<\cdots<k_m$ and that for each string $(i_1,i_2,\dots,i_m)$ we have established the following objects with the following properties:
\begin{enumerate}
\item A family of $(i_1,i_2,\dots,i_m)$-good $k_m$-rectangles. These families are disjoint for different strings $(i_1,i_2,\dots,i_m)$.
\item The closed \inv\ set $\hat S_{i_1,i_2,\dots,i_m}\subset\hat\X$, consisting of all arrays (allowed in $\hat\X$) whose first $k_m$ rows are concatenations of the $(i_1,i_2,\dots,i_m)$-good $k_m$-rectangles. These sets are nested: if $m\ge 2$ then $\hat S_{i_1,i_2,\dots,i_m}\subset\hat S_{i_1,i_2,\dots,i_{m-1}}$.
\item Isomorphic extensions
\begin{gather*}
\phi_{i_1,i_2,\dots,i_{m-1}}:\hat X_{i_1,i_2,\dots,i_{m-1}}\to X,\\ \theta_{i_1,i_2,\dots,i_m}:\hat X_{i_1,i_2,\dots,i_{m-1}}\to\hat X_{i_1,i_2,\dots,i_m},\\
\phi_{i_1,i_2,\dots,i_m}:\hat X_{i_1,i_2,\dots,i_m}\to X, \\
\text{such that \ } \phi_{i_1,i_2,\dots,i_m}\circ\theta_{i_1,i_2,\dots,i_m}=\phi_{i_1,i_2,\dots,i_{m-1}},
\end{gather*}
where $\hat X_{i_1,i_2,\dots,i_m}\subset\hat S_{i_1,i_2,\dots,i_m}$.
\item The associated map on measures
$$
\Phi_{i_1,i_2,\dots,i_m}=(\phi^*_{i_1,i_2,\dots,i_m})^{-1}:\K\to\hat\K_{i_1,i_2,\dots,i_m}=
\M^{\mathsf e}_\sigma(\hat X_{i_1,i_2,\dots,i_m}).
$$
%The sets $\hat\K_{i_1,i_2,\dots,i_m}$ are disjoint for different strings $(i_1,i_2,\dots,i_m)$.
\item The set $\hat K_{i_1,i_2,\dots,i_m}=\Phi_{i_1,i_2,\dots,i_{m-1}}(K_{i_1,i_2,\dots,i_m})$ of diameter less than $\eps_m$, such that $\M_\sigma(\hat S_{i_1,i_2,\dots,i_m})$ is contained in the $\eps_m$-neighborhood of the closed convex hull of $\hat K_{i_1,i_2,\dots,i_m}$. In particular, for any $\mu\in K_{i_1,i_2,\dots,i_m}$, we have
$$
d^*(\Phi_{i_1,i_2,\dots,i_{m-1}}(\mu),\Phi_{i_1,i_2,\dots,i_m}(\mu))< 2\eps_m.
$$
\item For each $\hat\mu\in\hat K_{i_1,i_2,\dots,i_m}$ the measure of points $\hat x$ whose central $k_m$-rectangle is changed by $\theta_{i_1,i_2,\dots,i_m}$ is less than $\eps_m$.
\end{enumerate}
\smallskip

Observe that the conditions (1)--(6) are fulfilled for $m=1$. Indeed, (1), (2), (4) and (6) are explicit. Further, in step $1$ the string $(i_1,i_2,\dots,i_{m-1})$ is empty, hence $\phi_{i_1,i_2,\dots,i_{m-1}}=\phi$ and $\Phi_{i_1,i_2,\dots,i_{m-1}}=\Phi$, which makes (3) and (5) consistent with step 1.

Note also that the disjointness in (1) implies that the sets $\hat S_{i_1,i_2,\dots,i_m}$ are disjoint for different strings $(i_1,i_2,\dots,i_m)$. As a consequence, also the extensions $\hat X_{i_1,i_2,\dots,i_m}$ are disjoint and so are the sets of measures $\hat\K_{i_1,i_2,\dots,i_m}$ (we will refer to this fact as condition (7)).
\smallskip

The construction of the analogous objects for $m\!+\!1$ consists in almost verbatim rewriting step 1. We establish a number $r_{m+1}$ such that each of the atoms $K_{i_1,i_2,\dots,i_m}$ of $\P_m$ can be partitioned into $r_{m+1}$ clopen subsets
$K_{i_1,i_2,\dots,i_{m+1}}$ (with $i_{m+1}\in\{1,2,\dots,r_{m+1}\}$), of diameters not exceeding $\frac{\diam(\K)}{m+1}$. We denote the resulting partition by $\P_{m+1}$.
We may arrange the atoms of $\P_{m+1}$ to be so small that for all possible strings $(i_1,i_2,\dots,i_{m+1})$ the sets
$\hat K_{i_1,i_2,\dots,i_{m+1}}=\Phi_{i_1,i_2,\dots,i_m}(K_{i_1,i_2,\dots,i_{m+1}})$ have diameters smaller than $\eps_{m+1}$. It is obvious that these sets are disjoint for different strings $(i_1,i_2,\dots,i_{m+1})$, and so are their closed convex hulls. We let $\gamma_{m+1}<\eps_{m+1}$ be smaller than half of the smallest distance between these hulls.

By Corollary~\ref{corol} and Proposition \ref{apcon}, there exists $k_{m+1}$ so large that
\begin{enumerate}
	\item[(a)] for each $\hat\mu\in\hat\K_{i_1,i_2,\dots,i_m}$ the measure of the set of points $\hat x$ whose central $k_{m+1}$-rectangle is $\gamma_{m+1}$-close to $\hat\mu$ exceeds $1-\gamma_{m+1}$,
	\item[(b)] there exists at least one $k_{m+1}$-rectangle of each length $l_{k_{m+1}}$ and $l_{k_{m+1}}\!+\!1$, which is $\gamma_{m+1}$-close to $\hat\mu$,
	\item[(c)] if $k_{m+1}$-rectangles $R_1,R_2,\dots,R_q$ satisfy $d^*(R_l,\hat\mu_l)<\gamma_{m+1}$, for some \im s $\hat\mu_l$, $l=1,2,\dots,q$, $q\in\N$, then $d^*(R,\bar\mu)<\eps_{m+1}$, where $R$ is the concatenation $R_1R_2\dots R_q$ and $\bar\mu$ is some convex combination of the measures $\hat\mu_l$.
\end{enumerate}

The further procedure will be described for a fixed string $(i_1,i_2,\dots,i_{m+1})$. Each of the $k_{m+1}$-rectangles appearing in $\hat X_{i_1,i_2,\dots,i_m}$, falling in the $\gamma_{m+1}$-neighborhood of $\hat K_{i_1,i_2,\dots,i_{m+1}}$ will be called \emph{$(i_1,i_2,\dots,i_{m+1})$-good}, all other $k_{m+1}$-rectangles are called \emph{$(i_1,i_2,\dots,i_{m+1})$-bad}. By the choice of $\gamma_{m+1}$, the families of
$(i_1,i_2,\dots,i_{m+1})$-good $k_{m+1}$-rectangles are disjoint for different strings
$(i_1,i_2,\dots,i_{m+1})$ (inductive condition (1)).
By (b), there exist $(i_1,i_2,\dots,i_{m+1})$-good $k_{m+1}$-rec\-tan\-gles of each length $l_{k_{m+1}}$ and $l_{k_{m+1}}\!+\!1$. We select one $(i_1,i_2,\dots,i_{m+1})$-good $k_{m+1}$-rectangle of each length and we call these two rectangles the \emph{$(i_1,i_2,\dots,i_{m+1})$-tabbed $k_{m+1}$-rectangles}. In every array $\hat x$ appearing in $\hat X_{i_1,i_2,\dots,i_m}$ we replace all $(i_1,i_2,\dots,i_{m+1})$-bad $k_{m+1}$-rectangles by the $(i_1,i_2,\dots,i_{m+1})$-tabbed $k_{m+1}$-rec\-tan\-gles of matching lengths. The modified system $\hat X_{i_1,i_2,\dots,i_m}$ will be denoted by $\hat X_{i_1,i_2,\dots,i_{m+1}}$, and the modification itself is $\theta_{i_1,i_2,\dots,i_{m+1}}: \hat X_{i_1,i_2,\dots,i_m}\to \hat X_{i_1,i_2,\dots,i_{m+1}}$. By (a), the set of arrays whose central $k_{m+1}$-rectangle is changed by $\theta_{i_1,i_2,\dots,i_{m+1}}$ has, for each $\hat\mu\in\hat K_{i_1,i_2,\dots,i_{m+1}}$, a measure smaller than~$\gamma_{m+1}$, which is less than $\eps_{m+1}$ (inductive condition (6)).

At the same time, we denote by $S_{i_1,i_2,\dots,i_{m+1}}$ the closed \inv\ subset of $\hat\X$ consisting of all possible arrays allowing any concatenations of $(i_1,i_2,\dots,i_{m+1})$-good $k_{m+1}$-rectangles in the first $k_{m+1}$-rows. Clearly, $\hat X_{i_1,i_2,\dots,i_{m+1}}\subset \hat S_{i_1,i_2,\dots,i_{m+1}}$ and the sets $\hat S_{i_1,i_2,\dots,i_{m+1}}$ are disjoint for different strings $(i_1,i_2,\dots,i_{m+1})$.
Recall that, by definition, all $(i_1,i_2,\dots,i_{m+1})$-good $k_{m+1}$-rectangles appear in $\hat X_{i_1,i_2,\dots,i_m}$. Thus, by congruency of the system of markers, every such $k_{m+1}$-rectangle has, in its first $k_m$ rows, a concatenation of $(i_1,i_2,\dots,i_m)$-good $k_m$-rectangles. This implies that  $\hat S_{i_1,i_2,\dots,i_{m+1}}\subset\hat S_{i_1,i_2,\dots,i_m}$ (inductive condition (2)).

As in step 1, the maps
\begin{gather*}
\phi_{i_1,i_2,\dots,i_m}:\hat X_{i_1,i_2,\dots,i_m}\to X,\\ \theta_{i_1,i_2,\dots,i_m}:\hat X_{i_1,i_2,\dots,i_m}\to\hat X_{i_1,i_2,\dots,i_{m+1}},\\
\phi_{i_1,i_2,\dots,i_{m+1}}:\hat X_{i_1,i_2,\dots,i_{m+1}}\to X,
\end{gather*}
where $\phi_{i_1,i_2,\dots,i_{m+1}}$ consists in erasing all markers and reconstructing the initial $k\!+\!1$ rows by amalgamations, satisfy $\phi_{i_1,i_2,\dots,i_{m+1}}\circ\theta_{i_1,i_2,\dots,i_{m+1}}=\phi_{i_1,i_2,\dots,i_m}$. This makes all these maps isomorphic extensions (inductive condition (3)). In particular, $\hat X_{i_1,i_2,\dots,i_{m+1}}$ is an isomorphic extension of $X$. The inverse of the adjoint map, $\Phi_{i_1,i_2,\dots,i_{m+1}}=(\phi_{i_1,i_2,\dots,i_{m+1}}^*)^{-1}:\K\to\hat\K_{i_1,i_2,\dots,i_{m+1}}=\M^{\mathsf e}_\sigma(\hat X_{i_1,i_2,\dots,i_{m+1}})$ is an affine homeomorphism sending measures to their isomorphic copies (inductive condition (4)). Since the sets $\hat\K_{i_1,i_2,\dots,i_{m+1}}$ with different strings $(i_1,i_2,\dots,i_{m+1})$ consist of measures supported by the disjoint systems $\hat S_{i_1,i_2,\dots,i_{m+1}}$, they are disjoint (condition (7)).

All elements of $\hat S_{i_1,i_2,\dots,i_{m+1}}$ have in the first $k_{m+1}$ rows concatenations of only $(i_1,i_2,\dots,i_{m+1})$-good $k_{m+1}$-rectangles, therefore, by (c), all \im s supported by $\hat S_{i_1,i_2,\dots,i_{m+1}}$ (in particular all elements of $\hat\K_{i_1,i_2,\dots,i_{m+1}}$) lie within the $\eps_{m+1}$-neighborhood of the closed convex hull of $K_{i_1,i_2,\dots,i_{m+1}}$. By convexity of $d^*$, this hull has diameter is smaller than $\eps_{m+1}$. As a consequence, for any $\mu\in K_{i_1,i_2,\dots,i_{m+1}}$, we have $d^*(\Phi_{i_1,i_2,\dots,i_m}(\mu),\Phi_{i_1,i_2,\dots,i_{m+1}}(\mu))< 2\eps_{m+1}$
(inductive condition (5)).
\smallskip

This concludes the induction. It remains to analyze what happens along the infinite paths $(i_1,i_2,i_3,\dots)\in\prod_{m=1}^\infty\{1,2,\dots,r_m\}$. First of all we observe, that
this infinite product is homeomorphic to $\K$ via the map that associates to each path
$(i_1,i_2,i_3,\dots)$ the unique measure $\mu_{i_1,i_2,i_3,\dots}$ determined by the
equality
$$
\bigcap_{m=1}^\infty K_{i_1,i_2,\dots,r_m} = \{\mu_{i_1,i_2,i_3,\dots}\}.
$$
Let us now fix a path $(i_1,i_2,i_3,\dots)$ and let us abbreviate $\mu_{i_1,i_2,i_3,\dots}$ as $\mu\in\K$ and $\Phi(\mu)$ as $\hat\mu\in\hat\K$. The property~(5) and summability of the \sq\ $(\eps_m)_{m\ge1}$ imply that the \sq\ of images $\Phi_{i_1,i_2,\dots,i_m}(\mu)$ converges. We let $\Psi(\mu)$ denote the limit measure.
Since for each $m\ge 1$ we have $\hat\mu\in\hat K_{i_1,i_2,\dots,i_m}$, the property~(6) implies that the passage $\theta_{i_1,i_2,\dots,i_{m+1}}$ from $\hat X_{i_1,i_2,\dots,i_m}$ to $\hat X_{i_1,i_2,\dots,i_{m+1}}$ changes symbols in the zero column in arrays $\hat x\in\hat X$ which constitute a set of measure $\hat\mu$ less than~$\eps_{m+1}$. By the Borel--Cantelli Lemma, the set of arrays $\hat x$ such that the zero column is changed during the entire inductive construction only finitely many times has measure $\hat\mu$ equal to 1. Obviously, the set $\hat X_{\hat \mu}\subset\hat X$ of arrays such that any column is changed only finitely many times also has measure 1. For any array $\hat x\in\hat X_{\hat\mu}$, the images
$\theta_{i_1,i_2,\dots,i_m}\circ\cdots\circ\theta_{i_1,i_2}\circ\theta_{i_1}(\hat x)\in\hat X_{i_1,i_2,\dots,i_m}$ converge (coordinatewise) to some array (with markers) which we denote by $\theta_\mu(\hat x)$. We also let $\tilde X_\mu=\theta_\mu(\hat X_\mu)$. The set
$\tilde X_\mu$ need not be closed but it supports the measure $\tilde\mu=\theta^*_\mu(\hat\mu)$.

Let $\hat x\in\hat X_{\mu}$. Since in each column only finitely many entries are changed by $\theta_\mu$, the original element $x$ can be reconstructed from $\theta_\mu(\hat x)$ by erasing all markers and applying suitable amalgamations. If we denote this reconstruction by $\phi_\mu$ then, on $\hat X_\mu$, we have $\phi_\mu\circ\theta_\mu=\phi$. This proves that all these maps are isomorphic extensions, and thus the measures $\mu$ and $\tilde\mu$ are isomorphic.

Consider an arbitrary rectangle $R$ and a $\hat\mu$-typical array $\hat x$. By the ergodic theorem, the frequency of the occurrence of $R$ in $\theta_\mu(\hat x)$ equals $\tilde\mu([R])$. On the other hand, it is seen (by summability of $(\eps_m)_{m\ge1}$) that it also equals the limit of the frequencies of occurrence of $R$ in the arrays $\theta_{i_1,i_2,\dots,i_m}\circ\cdots\circ\theta_{i_1,i_2}\circ\theta_{i_1}(\hat x)$, which are, in turn, equal to $\Phi_{i_1,i_2,\dots,i_m}(\mu)([R])$ and converge to $\Psi(\mu)([R])$. We conclude that the measures $\Psi(\mu)$ and $\tilde\mu$ coincide, and thus $\Psi(\mu)$ and $\mu$ are isomorphic.

\smallskip

Let us now analyze the structure of the sets $\hat S_{i_1,i_2,\dots,i_m}$. For each infinite path $(i_1,i_2,i_3,\dots)$ we define
$$
\hat S_{i_1,i_2,i_3,\dots} = \bigcap_{m=1}^\infty \hat S_{i_1,i_2,\dots,r_m}.
$$
By the ``nestedness'' stated in (2), this is a decreasing intersection, so it is nonempty, closed and \inv. The union
$$
\hat S_{(\infty)} = \bigcup_{(i_1,i_2,i_3,\dots)\in\prod_{m=1}^\infty\{1,2,\dots,r_m\}} \hat S_{i_1,i_2,i_3,\dots}
$$
is also an \inv\ set. Moreover, it is closed, because it equals the decreasing intersection
$$
\hat S_{(\infty)}=\bigcap_{m=1}^\infty\hat S_{(m)},
$$
where
$$
\hat S_{(m)}=\bigcup_{(i_1,i_2,\dots,i_m)\in\prod_{j=1}^m\{1,2,\dots,r_j\}} \hat S_{i_1,i_2,\dots,i_m}.
$$
In other words, we have constructed a \tl\ \ds\ $\hat S_{(\infty)}$ which splits as a union of subsystems $\hat S_{i_1,i_2,i_3,\dots}$ indexed by the infinite strings $(i_1,i_2,i_3,\dots)$, or, equivalently, by the measures $\mu\in\K$. Notice that the partition of $\hat S_{(\infty)}$ by the sets $\hat S_{i_1,i_2,i_3,\dots}$ is \usc.
This follows from the fact that the partition is a common refinement of the clopen (hence \usc) partitions by the sets $\hat S_{(\infty)}\cap\hat S_{i_1,i_2,\dots,i_m}$.\footnote{A common refinement of \usc\ partitions is \usc.}

In view of Theorem \ref{doweis}, the last thing we need to show in order to prove that $\hat S_{(\infty)}$ is the desired system $\hat X_{\mathsf{pu}}$ is, that each of the subsystems $\hat S_{i_1,i_2,i_3,\dots}$ is uniquely ergodic and supports a measure isomorphic to $\mu_{i_1,i_2,i_3,\dots}$.

Recall that by the property~(5), $\M_\sigma(\hat S_{i_1,i_2,\dots,i_m})$ is contained in the $\eps_m$-neigh\-bor\-hood of the closed convex hull of $\hat K_{i_1,i_2,\dots,i_m}$, which is a set of diameter less than $\eps_m$, thus the diameter of $\M_\sigma(\hat S_{i_1,i_2,\dots,i_m})$ does not exceed $3\eps_m$. This implies that $\M_\sigma(\hat S_{i_1,i_2,i_3,\dots})$ has diameter zero, i.e.\ it is a singleton. On the other hand, for $\mu=\mu_{i_1,i_2,i_3,\dots}$ we have, for each $m\ge 1$,
$$
\mu\in K_{i_1,i_2,\dots,i_m}, \text{ \ hence \ } \Phi_{i_1,i_2,\dots,i_{m-1}}(\mu)\in\hat K_{i_1,i_2,\dots,i_m}.
$$
Because $\M_\sigma(\hat S_{i_1,i_2,\dots,i_m})$ is contained in the $\eps_m$-neighborhood of
$K_{i_1,i_2,\dots,i_m}$, we conclude that the limit measure $\Psi(\mu)$ belongs to the intersection $\M_\sigma(\hat S_{i_1,i_2,i_3,\dots})$. Since $\hat S_{i_1,i_2,i_3,\dots}$ is uniquely ergodic, we have just identified its unique \im\ as $\Psi(\mu)$ which, as we have already proved, is isomorphic to $\mu$. The proof is now complete.
\end{proof}

\subsection{A strictly uniform model in the general case}

In this section we consider the general case, where the set $\K$ need not be \zd\ or compact, may have isolated points and contain periodic measures. We will construct a pure (with respect to $\K$) strictly uniform model of the initial system \xsmt. There is, however, one (trivial) exceptional case, when such a model (even a more general pure uniform model) cannot exist. The situation is captured in the precise formulation of our main theorem, which reads as follows:

\begin{thm}\label{main}
Let \xsmt\ be a \mpt\ and let $\K\subset\mtxe$ be a Borel-measurable set such that $\xi(\K)=1$, where $\xi$ is the ergodic decomposition of~$\mu$. A \zd, pure (with respect to $\K$), strictly uniform model $(X_{\mathsf{psu}},\sigma)$ of \xsmt\ exists if and only if it is {\bf not} the following ``\emph{exceptional case}'':
\begin{itemize}
	\item $\K$ consists of periodic measures only, and
	\item the set $\mathsf{Per}(\K)$ of the minimal periods of the measures in $\K$ contains
	an infinite \sq\ without a common divisor in $\Per(\K)$ (we agree that $1$ is a common
	divisor of all natural numbers).
\end{itemize}
\end{thm}

\begin{proof} The proof is lengthy and tedious, but modulo the application of Theorem~\ref{case1} it is a combination of relatively simple techniques. Let us begin by considering the ``exceptional case''. The limit of a weakly-star convergent \sq\ of periodic measures is either aperiodic or periodic with a period that divides all but finitely many periods in the \sq. So, in the exceptional case the set of \im s of any hypothetical pure model would contain a \sq\ without a limit. Since $\mtx$ is always compact, this cannot happen. This is why the ``exceptional case restriction'' cannot be relaxed.
\smallskip

We pass to the ``positive'' part of the proof. The \emph{general  scheme} is as follows: We replace $\K$ by its subset $\K'$ satisfying $\xi(\K')=1$, which is a countable union of some disjoint sets, one of which is a singleton: $\K'=\{\mu_0\}\cup\bigcup_{n\ge1}\K_n$. Then we create a strictly ergodic model $X_0$ for $\mu_0$ and strictly uniform systems $X_n\subset\X$ ($n\ge1$), pure with respect to the sets $\K_n$ (we will call them \emph{partial models}). Finally, we find disjoint conjugate copies $\bar X_n$ of the systems $X_n$, contained in the independent joining version of the universal system $\bar\X$, arranged in such a way that they converge (in the Hausdorff metric) to $X_0$. Then we define $X_{\mathsf{psu}}=\bigcup_{n\ge0}\bar X_n$. This set is closed, shift-invariant, and it is a disjopint union of strictly ergodic subsystems. By the necessity in Theorem \ref{doweis}, the partition into strictly ergodic subsystems is \usc\ individually on each partial model $\bar X_n$. Since the sets $\bar X_n$ converge to a strictly ergodic system, it is elementary to see that the partition of the entire set $X_{\mathsf{psu}}$ into strictly ergodic subsystems is \usc\ as well. By the sufficiency in Theorem \ref{doweis}, this implies that $X_{\mathsf{psu}}$ is strictly uniform. Thus, $X_{\mathsf{psu}}$ is a pure with respect to $\K$, strictly uniform model of \xsmt, as required.
\smallskip

With this scheme in mind, we first take care of periodic measures belonging to $\K$ (if they exist). They can be classified (up to isomorphism) by their minimal periods. For each period $p\in\Per(\K)$ the collection $\K^{\per}_p$ of all $p$-periodic measures in $\K$ is measurable and its cardinality $|\K^\per_p|$ is either finite or countable, or that of a continuum. In either case, there exists a strictly uniform (compact) \zd\ system $X^{\per}_p\subset\X$ consisting of $|\K^\per_p|$ $p$-periodic orbits. Since both $\K^\per_p$ and $\M^{\mathsf e}_\sigma(X^{\per}_p)$ are standard measurable spaces of equal cardinalities, by the Borel Isomorphism Theorem (see e.g. \cite{Ku}), there
exists a Borel-measurable bijection between them, and thus $X^{\per}_p$ is a strictly uniform system, pure with respect to $\K^\per_p$ (these will be our partial models).

\smallskip
We will now finish the proof in case $\K$ consists of periodic measures only, but it is not the ``exceptional case''. Then there exists a period $p_0\in\Per(\K)$ which divides all but finitely many periods in $\Per(\K)$. In this situation, we select one $p_0$-periodic measure $\mu_0$ and model it by a $p_0$-periodic orbit $X_0\subset\X$. We may need to redefine $X^{\Per}_{p_0}$ so that it is a model of $\K^{\Per}_{p_0}\setminus\{\mu_0\}$ (this is necessary only is $\K^{\Per}_{p_0}$ is finite). According to the general scheme, it now suffices to replace the strictly uniform models $X^\Per_p$ with $p$ such that $p_0|p$, by their conjugate copies contained in $\bar\X$ in such a way that they converge to $X_0$. This is done by a very simple manipulation of rows in the array representations of $X^\Per_p$. First, it is easy to see that for each $p$ there exists a clopen set visited by every orbit in $X^\Per_p$ precisely once in every $p$ iterates. This allows one to equip every array $x\in X^\Per_p$ with $p$-periodically distributed markers.
We can choose any row to place these markers, and we decide to place them in the $p$th row. Likewise, we equip the elements of $X_0$ with $p_0$-periodic markers in the $p_0$th row. Next, for each $p$ such that $p_0|p$ we modify all arrays $x\in X^\Per_p$ by replacing the initial $p$ rows by the initial $p$ rows of the unique $x_0\in X_0$ whose markers are congruent with those in $x$. Then we may remove the markers, so that the resulting arrays remain in $\bar\X$. The modification is deterministic, shift-equivariant, continuous, and since $X^\Per_p\subset\X$ and we alter only finitely many rows, the modification is invertible (the original initial rows can be reconstructed by amalgamations). Thus, we have created a conjugate copy $\bar X^\Per_p\subset\bar\X$ of $X^\Per_p$. For each $p$, $\bar X^\Per_p$ consist of $p$-periodic orbits, so, for different parameters $p$ the systems $\bar X^\Per_p$ are disjoint. It is clear, by the construction, that the systems $\bar X^\Per_p$ converge to $X_0$. This completes the construction of $X_{\mathsf{psu}}$ in the purely periodic (but not ``exceptional'') case.
\smallskip

Having dealt with the purely periodic case, from now on we will assume that $\K$ contains at least one aperiodic measure. We select one such measure and denote it by $\mu_0$. By a theorem of Lehrer (\cite{L}), there exists a zero-dimensional strictly ergodic and \tl ly mixing model $X_0\subset\X$ of $\mu_0$. The system $X_0$ will serve as the accumulation point for the rest of our model consisting of both periodic and aperiodic partial models.

So far we have created partial models for the measures in $\{\mu_0\}\cup\bigcup_{p\in\per(\K)}\K^\per_p$ and we still need to build partial models for the measures contained in
$$
\K^\triangle=\K\setminus\Bigl(\{\mu_0\}\cup\bigcup_{p\in\per(\K)}\K^\per_p\Bigr).
$$
According to the definition of a pure model, we are allowed to replace $\K^\triangle$ by its subset $\K^\bigtriangledown$ such that $\xi(\K^\bigtriangledown)=\xi(\K^\triangle)$. By regularity of the measure $\xi$ and elementary facts from topology, $\K^\triangle$ contains such a subset $\K^\bigtriangledown$ which is a disjoint at most countable union of closed \zd\ sets, $\K^\bigtriangledown=\bigcup_{s\in\mathbb S}\K^\bigtriangledown_s$ such that 
$\xi(\K^\bigtriangledown_s)>0$ for each $s\in\mathbb S$. By the Cantor--Bendixson Theorem, each set $\K^\bigtriangledown_s$ decomposes as a disjoint union of an at most countable set $\K^\omega_s$ and a perfect~set~$\K^\star_s$. Any zero-dimensional perfect set is a \tl\ Cantor set, and hence so~is~$\K^\star_s$.

We start by considering the countably many points of $\K^\omega=\bigcup_{s\in\mathbb S}\K^\omega_s$ (if they exist). We focus exclusively on those which are of positive measure $\xi$. By the Jewett--Krieger Theorem, each of them has a strictly ergodic \zd\ model contained in $\X$ and we can easily arrange these models to be different (hence disjoint) for different measures in $\K^\omega$ (even if some different measures in $\K^\omega$ are isomorphic), so that the map associating to each measure from $\K^\omega$ its model is a bijection. We denote these models by $X^\omega_j$, where $j$ ranges over some at most countable set $\mathbb J$.

Finally, we shall deal with the Cantor sets $\K^\star_s$. We redifine $\mathbb S$ to be the set of indices $s$ of only these sets $\K^\star_s$ which have positive measure $\xi$. By Theorem~\ref{case1}, for each $s\in\mathbb S$, there exist a \zd\ uniform system $X_{\!s}$ such that $\M^\mathsf e_\sigma(X_{\!s})$, from now on denoted by $\K_s$, and $\K_s^\star$ are indistinguishable. The conditional measure $\xi|\K^\star_s$ ($\xi$ restricted to subsets of $\K_s^\star$ and then normalized) corresponds via the indistinguishability to a probability measure supported by $\K_s$, which, by a slight abuse of notation, we will denote
by $\xi|\K_s$. Clearly, we can assume that the systems $X_{\!s}$ are subsystems of the universal system $\X$, and we can easily arrange them to be disjoint from each-other, as well as from the systems~$X^\omega_j$ (their disjointness from the periodic partial models $X^{\per}_p$ is automatic).

Our next goal is to replace each of the uniform systems $X_{\!s}$ by an at most countable disjoint union of strictly uniform ones. We achieve this using another regularity trick. It is well known (and easy to prove) that the mapping $\mu\mapsto\supp(\mu)$ assigning to Borel measures on a compact metric space $X$ their topological supports viewed as elements of $2^X$ (the collection of all compact subsets of $X$ equipped with the Hausdorff distance) is measurable (it is in fact lower semicontinuous). When applied to ergodic measures of the uniform system $X_{\!s}$ (where, by Theorem \ref{doweis}, each ergodic measure is supported by a separate uniquely ergodic subsystem), the map $\supp$ becomes injective, and hence bijective onto its image. By the Lusin-Souslin Theorem (see~\cite[Theorem 15.1 and Corollary 15.2]{K}, or \cite{Ku}), any Borel measurable bijection has a measurable inverse, which implies that the image $\mathsf S_s=\supp(\K_s)$ is a Borel subset of $2^{X_{\!s}}$. Let $\xi'_s$ denote the measure on $2^{X_{\!s}}$ which is the image of the measure $\xi|\K_s$ by the map adjoint to $\supp$. By regularity, there exists a subset of $\mathsf S'_s\subset\mathsf S_s$ such that $\xi'_s(\mathsf S'_s)=1$ and $\mathsf S'_s$ is a disjoint countable union of some closed subsets $\mathsf S_{s,t}$:
$$
\mathsf S'_s=\bigcup_{t\in\N}\mathsf S_{s,t}
$$
(for uniformity of our description, even if some sets $\mathsf S_s$ are closed, we replace them by a disjoint union of countably many closed subsets).
Let $\K_{s,t}\subset\K_s$ denote the preimage of $\mathsf S_{s,t}$ by the bijection $\supp$ and let
$\K^\star_{s,t}$ be the subset of $\K^\star_s$ corresponding to $\K_{s,t}$ via indistinguishability between $\K_s$ and $\K^*_s$. Finally, let
$$
\K^{\star\star}_s=\bigcup_{t\in\N}\K^\star_{s,t},
$$
where the union is obviously  disjoint. The map adjoint to $\supp$ sends measures on $\K_s$ to measures on $2^{X_{\!s}}$ and it clear that it sends $\xi|\K_s$ to a measure $\xi_s'$ satisfying
$$
\sum_{t\in\N}\xi'_s(\mathsf S_{s,t})=\xi'_s(\mathsf S'_s)=1.
$$
Thus, we have
$$
\xi(\K^{\star\star}_s)=\sum_{t\in\N}\xi(\K^\star_{s,t})=\xi(\K^\star_s).
$$

Let
$$
X_{s,t} =\bigcup\mathsf S_{s,t} = \bigcup_{\mu\in\K_{s,t}}\supp(\mu)\subset X_{\!s}.
$$
It is obvious that the set $X_{s,t}$ is invariant and closed (the union of a family of compact sets which is closed in the Hausdorff metric, is closed), and the set of ergodic measures carried by $X_{s,t}$ equals $\K_{s,t}$. Since, in a uniquely ergodic system, the support of the (unique) \im\ is minimal, the set $X_{s,t}$ is a union of minimal sets. Note also that the sets $X_{t,s}$ are disjoint for different indices $t$. As $X_{s,t}$ is contained in the uniform system $X_{\!s}$, it is uniform as well. In particular, $\K_{s,t}$ (and thus also $\K^\star_{s,t}$) is in fact compact (see \cite[Proposition~4.2(5) and Theorem~4.9]{DW}). We have just shown that $X_{s,t}$ is a strictly uniform system with $\M_\sigma^\mathsf e(X_{s,t})$ indistinguishable from $\K^\star_{s,t}$. In this manner, for each $s\in\mathbb S$, we have created a disjoint countable family of systems $X_{s,t}\subset X_{\!s}\subset\X$ ($t\in\N$) which are strictly uniform and pure with respect to the compact disjoint sets $\K^\star_{s,t}\subset\K^\star_s$ filling $\K^\star_s$ up to measure~$\xi$.

\smallskip
Up to now, we have constructed at most countably many pairwise disjoint strictly uniform partial models of three kinds: $X^\per_p$ ($p\in\Per(\K)$), $X^\omega_j$ ($j\in\mathbb J$), and $X_{s,t}$ ($s\in\mathbb S, t\in\N$). Our last step is finding conjugate copies of the partial models, which are disjoint and converge to the \tl ly mixing strictly ergodic system $X_0$. According to the argument explained in the \emph{general scheme}, this will end the proof. We will proceed assuming that each of the index sets $\Per(\K)$, $\mathbb J$ and $\mathbb S$ is infinite, otherwise we should just ignore the partial models corresponding to finite index sets. 
\smallskip

We begin by distributing some markers in the partial models. In each of them we will place the markers in only one (say, the first) row, and the procedure will be a conjugacy. And so, in each of the periodic systems $X^\per_p$ we have no choice but to put markers $p$-periodically (i.e.\ with only one gap size equal to $p$); the procedure was described a few paragraphs back. Next, we enumerate all countably many aperiodic systems of the form $X^\omega_j$ and $X_{s,t}$ in one \sq\ denoted $(X^\aper_n)_{n\in\N}$ and for each $n\in\N$ we place the markers with only two gap sizes, say $l_n$ and $l_n+1$; this procedure was also described earlier. We need to ensure that the \sq\ $l_n$ grows to infinity. Because the markers in $X^\per_p$ and $X^\aper_n$ are determined via a continuous and deterministic process, we can ``imagine'' them at the due places, without needing to actually put them in. This allows us to continue denoting the systems with markers by $X^\per_p$ and $X^\aper_n$, respectively.

It is now that we take the advantage of the \tl\ mixing property of $X_0$, which asserts that for any rectangle $B$ appearing in $X_0$ there exists $l(B)\in\N$ such that for any $l\ge l(B)$ the intersection $[B]\cap\sigma^l([B])\cap X_0$ is nonempty. For each $k\ge1$ we select a rectangle $B^{(k)}$ of dimensions $k\times 2k$ and we define
$$
l_k=k+\max\{l(B^{(i)}):i\le k\}.
$$
Clearly, the \sq\ $(l_k)_{k\ge 1}$ is strictly increasing, which enables us to
associate, to each $l\ge l_1$, the unique $k$ for which $l\in[l_k,l_{k+1}-1]$. We denote this unique $k$ by $k_l$. Further, to each $l\ge l_1$ we also associate the rectangle $B^{(k_l)}$, denote it by $B_l$ and call it the $l$th \emph{base rectangle} (the \sq\ of base rectangles $(B_l)_{l\ge l_1}$ may be constant on long intervals, but eventually both dimensions of $B_l$ grow to infinity). Observe that for each $l\ge l_1$, we have
$$
[B_l]\cap\sigma^l([B_l])\cap X_0\neq\emptyset \text{ \ and \ }[B_l]\cap\sigma^{l+1}([B_l])\cap X_0\neq\emptyset.
$$
These two facts imply that there exist some two elements $x^{(l)},\bar x^{(l)}\in X_0$ such that
$$
x^{(l)}_{[1,k_l]\times[0,2k_l-1]}=\bar x^{(l)}_{[1,k_l]\times[0,2k_l-1]}=
x^{(l)}_{[1,k_l]\times[l,l+2k_l-1]}=\bar x^{(l)}_{[1,k_l]\times[l+1,l+2k_l]}=B_l.
$$
We assign
$$
R_l = x^{(l)}_{[1,k_l]\times[k_l,l+k_l-1]} \text{ \ \ and \ \ }\bar R_l = \bar x^{(l)}_{[1,k_l]\times[k_l,l+k_l]}
$$
and call them the \emph{tabbed rectangles}. The dimensions of the tabbed rectangles are $k_l\times l$ and $k_l\times(l\!+\!1)$, respectively (see Figure \ref{fig2}).

\begin{figure}[H]
\includegraphics[width=10cm]{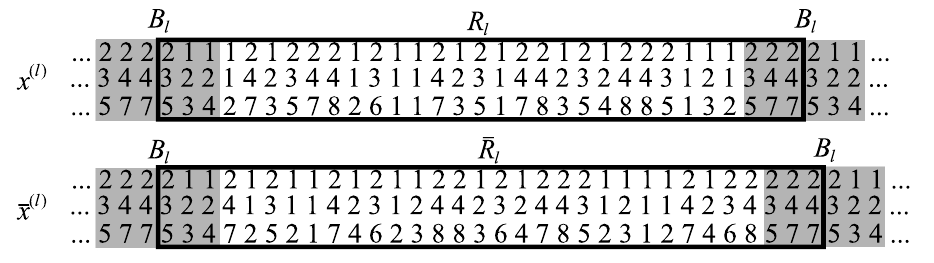}
\vspace{10pt}
\caption{The tabbed rectangles $R_l$ and $\bar R_l$ are shown in the black frames. Their lengths are $l$ and $l+1$, respectively. They start and end in the middle of copies of the base rectangle $B_l$ (shown in gray).}
\label{fig2}
\end{figure}

By the obvious inequality $l\ge k_l$, each of the tabbed rectangles starts with the complete right half of the base rectangle $B_l$ and ends with the complete left half of the same base rectangle. This implies that in any concatenation of two tabbed rectangles (i.e. in $R_lR_l$, $\bar R_lR_l$, $R_l\bar R_l$ or $\bar R_l\bar R_l$), at the junction there occurs the complete base rectangle $B_l$ positioned symmetrically around the junction. Because both the tabbed rectangles and the base rectangle appear in $X_0$, any sub-rectangle of the concatenation, of dimensions $k_l\times k_l$ (or smaller) also appears in $X_0$. We will refer to this property of the tabbed rectangles as being ``seamlessly stitchable''~(in~$X_0$).

We are ready to build the desired conjugate copies of the systems $X^\per_p$ and $X^\aper_n$. We start by modifying the periodic ones. For each $p\in\Per(\K)$, in each array $x\in X^\per_p$, we replace every rectangle of dimensions $k_p\times p$ (i.e. $k_l\times l$, where $l=p$), appearing in the initial $k_p$ rows between two consecutive markers (which are always $p$ positions apart), by the tabbed rectangle $R_p$ (the rectangle $\bar R_p$ will not be used in handling the periodic partial models).
Before we discuss the properties of the above modification of $X^\per_p$, we pass to modifying the aperiodic systems $X^\aper_n$. Recall that, for each $n\ge1$, the markers in $X^\aper_n$ appear with only two gaps, $l_n$ and $l_n+1$. In each array $x\in X^\aper_n$, we replace every rectangle of dimensions $k_{l_n}\times l_n$ or $k_{l_n}\times (l_n\!+\!1)$, appearing in the initial $k_{l_n}$ rows between two consecutive markers, by the tabbed rectangle $R_{l_n}$ or $\bar R_{l_n}$, according to the length. This concludes the construction.

The above modifications of the systems $X^\per_p$ and $X^\aper_n$ are deterministic, shift-equivariant, continuous and can be inverted using the amalgamations. So, the resulting systems, denoted by $\bar X^\per_p$ and $\bar X^\aper_n$, are conjugate to their respective counterparts (the systems $\bar X^\per_p$ and $\bar X^\aper_n$ are no longer contained in $\X$, nevertheless, they are contained in~$\bar\X$).

It remains to prove that both \sq s of systems $(\bar X^\per_p)_{p\in\Per(\K)}$ and $(\bar X^\aper_n)_{n\in\N}$ converge to $X_0$ (in the Hausdorff metric). By minimality of $X_0$,
it suffices to show that $X_0$ contains the upper limit sets of the considered \sq s (because then both the upper and lower limit sets are closed invariant subsets of $X_0$, so they must equal $X_0$). Consider any rectangle $R$ which appears in infinitely many of the systems $\bar X^\per_p$ (respectively, $\bar X^\aper_n$). What we need to show is that $R$ appears in~$X_0$.
For $p$ (respectively, $n$) large enough, both dimensions of $R$ are smaller than or equal to~$k_p$ (respectively,~$k_{l_n}$). Then $R$ is part of a rectangle of dimensions $k_p\times k_p$ (respectively $k_{l_n}\times k_{l_n}$) appearing in $\bar X^\per_p$ (respectively, $\bar X^\aper_n$) in the initial~$k_p$ (respectively~$k_{l_n}$) rows filled entirely by the tabbed rectangles. The ``seamless stitchability'' of the tabbed rectangles guarantees that no matter whether $R$ appears inside one tabbed rectangle or over a junction, it appears also in $X_0$. We have just proved the last ingredient needed to complete the proof of Theorem \ref{main}.
\end{proof}

We remark that without the base rectangles, which enable us to ``seamlessly stitch'' the tabbed rectangles, we could only get the systems $\bar X^\per_p$ and $\bar X^\aper_n$ to converge to some uniquely ergodic system containing $X_0$, but not exactly to $X_0$. Some rectangles appearing over the junctions of the tabbed rectangles might not appear in $X_0$ and ``survive'' in the limit system (although they would have measure zero there). This is why some version of \tl\ mixing seems necessary and Lehrer's result comes in perfectly handily.

\end{document}